\documentclass[12pt]{article}

\begin{document}
\begin{center}
{\large {\bf Chains of twists for symplectic Lie algebras}}\\[5mm]

{\sc {\bf David N. Ananikian}
\\ Theoretical
Department, St. Petersburg State University,\\ 198904, St. Petersburg,
Russia \\ \vskip0.5cm }

{\sc {\bf Petr P. Kulish}
\\ St. Petersburg Department of the Steklov Mathematical
\\ Institute, 191011,
St.Petersburg, Russia \\ \vskip0.5cm }

{\sc {\bf Vladimir D. Lyakhovsky}
\\ Theoretical
Department, St. Petersburg State University,\\ 198904, St. Petersburg,
Russia \\ \vskip0.5cm }

\end{center}

\vskip0.5cm

\begin{abstract}
Serious difficulties arise in the construction of chains of twists for
symplectic simple Lie algebras.
Applying the canonical chains of extended twists to deform the
Hopf algebras $U(sp(N))$ one is forced to deal only with improper chains
(induced by the $U(sl(N))$ subalgebras). In the present paper this problem
is solved. For chains of regular injections $U(sp(1))\subset ...\subset
U(sp(N-1))\subset U(sp(N))$ the sets of maximal extended jordanian twists $%
{\cal F}_{{\cal E}_k}$ are considered. We prove that there exists for $U(sp(N))$
the twist ${\cal F}_{{\cal B}_{0 \prec k}}$ composed 
of the factors ${\cal F}_{{\cal E}_k}$. 
It is demonstrated that the twisting procedure deforms the space of the primitive
subalgebra $sp(N-1)$. The recursive algorithm for such deformation is found.
This construction generalizes the results obtained for orthogonal classical
Lie algebras and demonstrates the universality of primitivization effect for
regular chains of subalgebras. For the chain of the maximal length the
twists ${\cal F}_{{\cal B}_{0\prec k^{\max }}\ }$become full, their
carriers contain
the
Borel subalgebra $B^{+}\left( sp(N)\right) $. Using such twisting procedures
one can obtain the explicit quantizations for a wide class of classical $r$%
-matrices. As an example the full chains of extended twists for $U(sp(3))$
is constructed.
\end{abstract}

\section{Introduction}

Quantizations of triangular Lie bialgebras with antisymmetric classical $r$%
-matrices $r=-r_{21}$ are defined by the twisting elements ${\cal F}=\sum
f_{\left( 1\right) }\otimes f_{\left( 2\right) }\in {\cal A}\otimes {\cal A}$
that satisfy the twist equations ~\cite{DR83}:
\begin{equation}
\label{twist-eq}
\begin{array}{l}
\left( {\cal F}\right) _{12}\left( \Delta \otimes {\rm id}%
\right) {\cal F}=\left( {\cal F}\right) _{23}\left( {\rm id}\otimes \Delta
\right) {\cal F}, \\ \left( \epsilon \otimes {\rm id}\right) {\cal F}=\left(
{\rm id}\otimes \epsilon \right) {\cal F}=1.
\end{array}
\end{equation}

In applications the knowledge of the twisting elements is quite important
giving twisted ${\cal R}$-matrices ${\cal R_F}={\cal F}_{21}{\cal R}{\cal F%
}^{-1}$ and coproducts
$$
\Delta _{{\cal F}}={\cal F}\Delta {\cal F}^{-1}.
$$

Only a few classes of twists can be written explicitly in the closed form ~%
\cite{R,OGI,GIA,KLM}. The most interesting of these twists are the so called
{\sl extended jordanian twists} ${\cal F}_{{\cal E}}$ ~\cite{KLM} that can
be defined on a special type of carrier algebras ${\bf L}$. These algebras
are the analogs of the enlarged Heisenberg algebra and can be found in any
simple Lie algebra $g$ of rank greater than 1. The minimal algebra of this
type is 4-dimensional:
$$
\begin{array}{l}
\begin{array}{l}
\left[ H,E\right] =E, \\
\left[ A,B\right] =E,
\end{array}
\quad
\begin{array}{l}
\left[ H,A\right] =\alpha A, \\
\left[ H,B\right] =\beta B,
\end{array}
\\
\qquad \qquad \alpha +\beta =1.
\end{array}
$$

The explicit form of the twisting factors in the extended jordanian twist $%
{\cal F}_{{\cal E}}$ is:
$$
\begin{array}{l}
\Phi _{
{\cal J}}=e^{H\otimes \sigma }, \\ \Phi _{
{\cal E}}=e^{A\otimes Be^{-\beta \sigma }}, \\
\end{array}
\qquad \sigma =\ln \left( 1+E\right).
$$

The first factor $\Phi _{{\cal J}}$ is the jordanian twist \cite{OGI}
corresponding to the classical $r$-matrix $r_{{\cal J}}=H\wedge E$. The
second one, $\Phi _{{\cal E}}$, satisfies the twist equations~\ref{twist-eq}
only on the deformed quasiprimitive costructures. In particular the
necessary coalgebra can be obtained performing the jordanian twisting by $%
\Phi _{{\cal J}}$ in the initial Heisenberg subalgebra:
$$
\begin{array}{l}
\Delta _{
{\cal J}}\left( A\right) =A\otimes e^{\alpha \sigma }+1\otimes A, \\ \Delta_{
{\cal J}}\left( B\right) =B\otimes e^{\beta \sigma }+1\otimes B, \\ \Delta _{%
{\cal J}}\left( E\right) =E\otimes e^\sigma +1\otimes E.
\end{array}
$$

The composition of the twists $\Phi _{{\cal J}}$ and $\Phi _{{\cal E}}$
defines the extended jordanian twist \cite{KLM}:
$$
{\cal F}_{{\cal E}}=e^{A\otimes Be^{-\beta \sigma }}e^{H\otimes \sigma }.
$$
Here the corresponding classical $r$-matrix is 
$r_{{\cal EJ}}=H\wedge E+A\wedge B$.

It was demonstrated in \cite{KLO} that the extended twists ${\cal F}_{{\cal E%
}}$ can be composed into {\em chains}. The latter are based on the sequences
of regular injections
\begin{equation}
\label{inject} g_p\subset g_{p-1}\ldots \subset g_1\subset g_0=g.
\end{equation}
To compose the sequence one must choose the {\sl initial} root
$\lambda^k_0$ in each root system $\Lambda \left( g_k \right) $.
Let $V_{\lambda_0^k}^{\perp }$ be the subspace orthogonal
to the initial root $\lambda^k_0$
in the system $\Lambda \left( g_k \right) $. The subalgebra $g_{k+1}$
is defined by the subsystem $\Lambda \left( g_{k+1}\right) =$ $\Lambda
\left( g_k \right) \cap V_{\lambda^k_0}^{\perp }$. Consider the set $\pi_k $
of the constituent roots for $\lambda^k_0$
\begin{equation}
\label{c-roots}
\begin{array}{l}
\pi_k =\left\{ \lambda ^{\prime },\lambda ^{\prime \prime }\mid \lambda
^{\prime }+\lambda ^{\prime \prime }=\lambda^k_0;\quad \lambda ^{\prime
}+\lambda^k_0,\lambda ^{\prime \prime }+\lambda^k_0\notin \Lambda \left(
g_k\right) \right\} \\
\pi_k =\pi_k ^{\prime }\cup \pi_k ^{\prime \prime };\qquad \pi_k ^{\prime }
=\left\{ \lambda ^{\prime }\right\} ,\pi_k^{\prime \prime }=\left\{ \lambda
^{\prime \prime }\right\} .
\end{array}
\end{equation}
It was shown that for the classical Lie algebras one can always find in $%
g_{\lambda^k_0}^{\perp }$ a subalgebra $g_{k+1}\subseteq
g_{\lambda^k_0}^{\perp }\subset g_k$ whose generators become primitive after
the twist ${\cal F}_{{\cal E}_k}$.

Such a process of primitivization of $g_{k+1}\subset g_{k}$ (called the {\em %
matreshka} effect) provides the possibility to compose chains of extended
twists of the type
\begin{equation}
\label{chain-ini}
\begin{array}{l}
{\cal F}_{{\cal B}_{0\prec p}}=\prod_{\lambda ^{\prime }\in \pi _p^{\prime
}}\exp \left\{ E_{\lambda ^{\prime }}\otimes E_{\lambda _0^p-\lambda
^{\prime }}e^{-\frac 12\sigma _{\lambda _0^p}}\right\} \cdot \exp
\{H_{\lambda _0^p}\otimes \sigma _{\lambda _0^p}\}\,\cdot \\ \prod_{\lambda
^{\prime }\in \pi _{p-1}^{\prime }}\exp \left\{ E_{\lambda ^{\prime
}}\otimes E_{\lambda _0^{p-1}-\lambda ^{\prime }}e^{-\frac 12\sigma
_{\lambda _0^{p-1}}}\right\} \cdot \exp \{H_{\lambda _0^{p-1}}\otimes \sigma
_{\lambda _0^{p-1}}\}\,\,\cdot \\
\ldots \\
\prod_{\lambda ^{\prime }\in \pi _0^{\prime }}\exp \left\{ E_{\lambda
^{\prime }}\otimes E_{\lambda _0^0-\lambda ^{\prime }}e^{-\frac 12\sigma
_{\lambda _0^0}}\right\} \cdot \exp \{H_{\lambda _0^0}\otimes \sigma
_{\lambda _0^0}\}\,.
\end{array}
\end{equation}

Some peculiarities were found in the case of orthogonal simple algebras ~%
\cite{Full, LSK}. Here the subalgebra $g_{\lambda _0^k}^{\perp }$ is
isomorphic to the direct sum $g_{\lambda _0^k}^{\perp }\approx sl(2)\oplus
g_{k+1}$. After being twisted by ${\cal F}_{{\cal E}_k}$ the costructure of
the first summand $sl(2)$ is nontrivially deformed while the second $%
g_{k+1}=so(N-4(k+1))$ remain primitive. The primitive summands form the
carrier space for the canonical chains \cite{KLO}. These chains are based on
the set of injections $so(N-4(k+1))\in so(N-4k)$. In \cite{Full} it was
proved that each algebra $U_{{\cal E}_k}(g_{\lambda _0^k}^{\perp })$ contains not
only the deformed subalgebra $U_{{\cal E}_k}(sl^k(2))$ but also the primitive Hopf
subalgebra $U(sl_G^k(2))$. Due to the fact that the constituent roots (\ref
{c-roots}) form the weight diagram for the representation of $so(N-4(k+1))$
in $g_{\lambda _0^k}^{\perp }$ the invariants of this representation can be
used to construct the generators of $sl_G^k(2)$. Thus the universal
character of the primitivization effect was confirmed. It was demonstrated
that chains of twists for orthogonal algebras can have the properties
similar to those of linear simple algebras: after $k$ subsequent steps one
can find the primitive subalgebra equivalent to $g_{\lambda
_0^k}^{\perp }$ though it is realized on a {\sl deformed carrier subspace}.

For the symplectic simple Lie algebras the situation is more complicated.
The coproducts of generators which correspond to the subalgebra $g_{\lambda
_0^k}^{\perp }\subset g_k$ are nontrivially deformed by the twist ${\cal F}_{%
{\cal E}_{k-1}{{\cal J}_{k-1}}}$. Hence a proper chain of twists can not
be constructed in a canonical way. Nevertheless as we shall demonstrate
below there exist
in $U_{{\cal F}_{{\cal E}_{k-1}{{\cal J}_{k-1}}}}\left( g_{k-1}\right) $
the primitive subalgebras equivalent to $g_{\lambda
_0^k}^{\perp }$. The equivalence map is realized by a specific nonlinear
transformation of basis. To clarify the algorithm we consider the particular
case of $U\left( sp(3)\right) $. This is the simplest example where the
specific structure of symplectic algebras appears.

\section{Chains of twists for $sp(3)$}

Consider the root system $\Lambda \left( sp(3)\right) $:
$$
\Lambda =\{\pm e_i\pm e_j,\quad \pm 2e_i\}\ (i,j=1,2,3).
$$
For the initial root $\lambda _0=2e_1$ the constituent roots are $\lambda
^{\prime }=e_1-e_i$ and $\lambda ^{\prime \prime }=e_1+e_i$. We shall use
the basis $\left\{H_{ii}, E_{\lambda_j}, F_{\lambda_j}| \lambda_j \in
\Lambda^+ \left( sp(3)\right) \right\}$. For the nonzero root $\lambda_j$
let us denote the generators $L_{\lambda_j}=L_{\pm e_i\pm e_j}$ by $L_{\pm i
\pm j}$. In this basis the Borel subalgebra $B^+(sp(N))$ is defined by the
relations
$$
\begin{array}{l}
\,[H_{ii},E_{n+n}]=\delta _{in}E_{n+n}, \\
\,[H_{ii},E_{n+m}]=\frac 12(\delta _{in}+\delta _{im})E_{n+m}, \\
\,[H_{ii},E_{n-m}]=\frac 12(\delta _{in}-\delta _{im})E_{n-m}, \\
\end{array}
\begin{array}{l}
\,[E_{i-j},E_{n+n}]=2\delta _{jn}E_{i+n}, \\
\,[E_{i-j},E_{n-m}]=\delta _{jn}E_{i-m}-\delta _{mi}E_{n-j}, \\
\,[E_{i-j},E_{n+m}]=\delta _{jn}E_{i+m}+\delta _{mi}E_{n+j}.
\end{array}
$$
The other $sp(N)$-commutators can be obtained using the Chevalley involution
$H_{ii}\rightarrow -H_{ii},\ E_{i\pm j}\rightarrow F_{i\pm j},\ F_{i\pm
j}\rightarrow E_{i\pm j}$.

Here the set of regular subalgebras 
$$
U(sp(1))\subset \dots \subset U(sp(k))\subset \dots \subset U(sp(N))
$$
coincides with the set of subalgebras $g_{\lambda^k_0}^{\perp }$
(compare with (\ref{inject})). Moreover
for the corresponding injections of the root systems
$$
\Lambda(sp(1))\subset \dots \subset \Lambda(sp(k))\subset \dots \subset
\Lambda(sp(N))
$$
the following property is true: the roots of $\Lambda(sp(k))$ are orthogonal
to any long root $\lambda_{0}\in \Lambda(sp(N))\backslash \Lambda (sp(k))$.
Notice that the chain of regular injections can be based on an arbitrary
long root in
$\Lambda \left( sp(N)\right)$. In our particular case the corresponding sets
are
$$
\begin{array}{c}
U(sp(1))\subset U(sp(2))\subset U(sp(3)), \\
\Lambda(sp(1))\subset \Lambda(sp(2)) \subset \Lambda(sp(3)).
\end{array}
$$
The property mentioned above (for the appropriate ordering of roots) is
quite simple: $\Lambda(sp(2))\,\bot\, 2e_{1}$ and $\Lambda(sp(1))\,\bot\, 
2e_{1}, 2e_{2}$.

\subsection{The first step -- the full extended twist}

We start the construction of the full chain of twists for $U\left(
sp(3)\right)$ by performing the jordanian twist with the carrier subalgebra
generated by $\left\{H_{11}, E_{1+1} \right\}$:
$$
\Phi _{J_1}=exp\{H_{11}\otimes \sigma _{1+1}\},\qquad \sigma
=\ln(1+E_{1+1}).
$$
The sets $\{2e_1,e_{1-2},e_{1+2}\}$, $\{2e_1,e_{1-3},e_{1+3}\}$ define two
extensions ${\cal E}^{\prime}$ and ${\cal E}^{\prime\prime}$ for $\Phi _{J_1}
$. So the full extended jordanian twist has the twisting element
\begin{equation}
\label{exttw} {\cal F}_{{\cal E}_1}=\underbrace{e^{E_{1-3} \otimes
E_{1+3}e^{-\frac 12\sigma _{1+1}}}}_{{\cal E}^{\prime\prime}}\underbrace{%
e^{E_{1-2}\otimes E_{1+2}e^{-\frac 12\sigma _{1+1}}}}_{{\cal E}%
^{\prime}}\cdot \underbrace{e^{H_{11}\otimes \sigma _{1+1}}}_J.
\end{equation}
Let us write down the corresponding twisted costructure of the subalgebra $%
g_{\lambda^0_0}^{\perp }=g_{(2e_1)}^{\perp }= U_{{\cal E}_{1}}(sp(2))$,
\begin{equation}
\label{costr}
\begin{array}{lcl}
\Delta_{{\cal E}_1} (E_{2+2}) & = & E_{2+2} \otimes 1 + 1\otimes E_{2+2}
+2E_{1+2}\otimes E_{1+2} e^{-
\frac{1}{2}\sigma_{1+1}}+ \\  &  & + E_{1+1}\otimes E^{2}_{1+2}
e^{-\sigma_{1+1}}, \\
\Delta_{{\cal E}_1} (E_{2+3}) & = & E_{2+3} \otimes 1 + 1\otimes E_{2+3}+
E_{1+2} \otimes E_{1+3} e^{-
\frac{1}{2}\sigma_{1+1}}+ \\  &  & + E_{1+3} \otimes E_{1+2} e^{-
\frac{1}{2}\sigma_{1+1}}+ E_{1+1} \otimes E_{1+2}E_{1+3} e^{-\sigma_{1+1}},
\\ \Delta_{{\cal E}_1} (E_{2-3}) & = & E_{2-3} \otimes 1 + 1\otimes E_{2-3},
\\
\Delta_{{\cal E}_1} (E_{3+3}) & = & E_{3+3} \otimes 1 + 1\otimes E_{3+3}
+2E_{1+3}\otimes E_{1+3} e^{-
\frac{1}{2}\sigma_{1+1}}+ \\  &  & + E_{1+1}\otimes E^{2}_{1+3}
e^{-\sigma_{1+1}}, \\
\Delta_{{\cal E}_1} (F_{2+2}) & = & F_{2+2} \otimes 1 + 1\otimes F_{2+2}
+2E_{1-2}\otimes E_{1-2} e^{-
\frac{1}{2}\sigma_{1+1}}- \\  &  & - E_{1+2}^{2}\otimes E_{1+1}
e^{-\sigma_{1+1}}, \\
\Delta_{{\cal E}_1} (F_{2+3}) & = & F_{2+3} \otimes 1 + 1\otimes F_{2+3} +
E_{1-2} \otimes E_{1-3} e^{-
\frac{1}{2}\sigma_{1+1}}+ \\  &  & + E_{1-3} \otimes E_{1-2} e^{-
\frac{1}{2}\sigma_{1+1}} - E_{1-2}E_{1-3}\otimes E_{1+1} e^{-\sigma_{1+1}},
\\ \Delta_{{\cal E}_1} (F_{2-3}) & = & F_{2-3} \otimes 1 + 1\otimes F_{2-3},
\\
\Delta_{{\cal E}_1} (F_{3+3}) & = & F_{3+3} \otimes 1 + 1\otimes F_{3+3}
+2E_{1-3}\otimes E_{1-3} e^{-
\frac{1}{2}\sigma_{1+1}}- \\  &  & - E_{1+3}^{2}\otimes E_{1+1}
e^{-\sigma_{1+1}}.
\end{array}
\end{equation}

To continue the construction of the chain we need a primitive subalgebra.
The generators $E_{2-3},F_{2-3}$ (and the corresponding Cartan element) are
primitive but these are not sufficient for a proper full chain \cite{KLO}.
It was demonstrated in \cite{KL} that the {\sl deformed carrier
subspaces},
algebraically equivalent to 
$U(g_{\lambda _0}^{\perp })$, can be found in $U_{{\cal E}}$. 
Applying the algorithms elaborated in \cite{KL} to the Hopf
subalgebra $U_{{\cal E}}(sp(2))$ with the costructure (\ref{costr}) we find
the following basis transformation.
$$
\begin{array}{lcl}
E^{\prime}_{i+i} & = & E_{i+i} - E^{2}_{1+i} e^{-\sigma_{1+1}}, \\
F^{\prime}_{i+i} & = & F_{i+i} - E^{2}_{1-i}, \\
E^{\prime}_{2+3} & = & E_{2+3} - E_{1+2}E_{1+3} e^{-\sigma_{1+1}}, \\
F^{\prime}_{2+3} & = & F_{2+3} - E_{1-2}E_{1-3}, \\
H^{\prime}_{ii} & = & H_{ii}.
\end{array}
\qquad i=2,3
$$
It can be checked explicitly that these new basic elements are primitive
with respect to $\Delta_{{\cal E}_1}$. They generate the subalgebra $%
sp^{\prime}(2)$.

\subsection{The next step -- the second extended twist}

Now we can temporary forget about the twist (\ref{exttw}). Consider the 
algebra $U(sp^{\prime}(2))$ 
on the deformed carrier space generated by $E^{\prime}$, $%
F^{\prime}$ and $H^{\prime}$. Introduce for it an independent root
system $\Lambda(sp^{\prime}(2)) $. In what follows we shall
drop the primes in the notation for the corresponding roots.

Let us perform in $sp^{\prime}(2) $ the jordanian twist based on the long
root $\lambda_{2+2}=2e_{2}\in \Lambda(sp^{\prime}(2)) $
\begin{equation}
\label{2-jord}
\begin{array}{c}
\Phi_{J_{2}}=\exp\{H_{22}\otimes\sigma^{\prime}_{2+2}\}, \\
\sigma^{\prime}_{2+2}= \ln (1+E^{\prime}_{2+2}).
\end{array}
\end{equation}

The set $\left\{ 2e_{2+2},\ e_{2-3},\ e_{2+3} \right\}$ indicates that there
exists the extension ${\cal E}_2$ for the twist (\ref{2-jord}):
$$
{\cal F}_{{\cal E}_2}=\underbrace {e^{E^{\prime}_{2-3}\otimes
E^{\prime}_{2+3} e^{-
\frac{1}{2}\sigma^{\prime}_{2+2}}}}_{{\cal E}} e^{H_{22}\otimes
\sigma^{\prime}_{2+2}}
$$

Applying this twist to $U(sp^{\prime}(2)) $ we see that
the subalgebra $sp^{\prime}(1) \in
sp^{\prime}(2)$ becomes nontrivially deformed. Again, the deformed carrier
space can be found. The following basis transformation gives the appropriate
primitive subalgebra $sp^{\prime \prime}(1) \in U_{{\cal E}_2{\cal E}%
_1}(sp(3))$:
$$
\begin{array}{lcl}
E^{\prime\prime}_{3+3} & = & E^{\prime}_{3+3} - (E^{\prime}_{2+3})^{2}
e^{-\sigma^{\prime}_{2+2}}, \\
F^{\prime\prime}_{3+3} & = & F^{\prime}_{3+3} - (E^{\prime}_{2-3})^{2}, \\
H^{\prime\prime}_{33} & = & H_{33}. \\
&  &
\end{array}
$$

\subsection{The last step -- the jordanian twist}

Now consider the subalgebra  $U(sp^{\prime\prime}(1)) $ on the
twice nontrivially
deformed space generated by $E^{\prime\prime}$, $F^{\prime\prime}$ and $%
H^{\prime\prime} $ and make the last step in our chain construction.

Perform in $U(sp^{\prime\prime}(1)) $ the jordanian twist based on the last
long root $\lambda_{3+3}=2e_{3} $,
$$
\Phi_{J_{3}}=\exp\{ H_{33}\otimes \sigma^{\prime\prime}_{3+3}\},
$$
$$
\sigma^{\prime\prime}_{3+3}= \ln(1+E^{\prime \prime}_{3+3}).
$$

There are no constituent roots for $\lambda_{3+3}$. So the last factor of
the chain is purely jordanian:
$$
{\cal F}_{{\cal E}_3}=\Phi_{J_{3}}.
$$

Notice, that in terms of the initial basic elements the deformed
generator, $%
E^{\prime\prime}_{3+3} $, looks like:
$$
E^{\prime\prime}_{3+3}=E_{3+3}-E_{1+3}^{2}
e^{-\sigma_{1+1}}-(E_{2+3}-E_{1+2}E_{1+3} e^{-\sigma_{1+1}})^{2}
e^{-\sigma^{\prime}_{2+2}}.
$$
Even for rather a small chain it appears extremely difficult to find such
complicated deformations of the carrier space without the recursive
procedure described above.

In this way the step by step transformations allow us to construct the
following chain of twists for $U(sp(3))$:
\begin{equation}
\label{f-chain}
\begin{array}{lcl}
{\cal F}_{{\cal G}_{0 \prec 2 }} & = & e^{H_{33}\otimes
\sigma^{\prime\prime}_{3+3}}\cdot \\
&  & \cdot e^{E^{\prime}_{2-3}\otimes E^{\prime}_{2+3} e^{-
\frac{1}{2}\sigma^{\prime}_{2+2}}} e^{H_{22}\otimes
\sigma^{\prime}_{2+2}}\cdot \\  &  & \cdot e^{E_{1-3}\otimes
E_{1+3}e^{-\frac 12\sigma _{1+1}}} e^{E_{1-2}\otimes E_{1+2}e^{-\frac
12\sigma_{1+1}}} e^{H_{11}\otimes \sigma _{1+1}}. \\
&  &
\end{array}
\end{equation}
Actually the set of injections was also deformed:
$$
U_{{\cal F}_{{\cal E}_3}}(sp(1)) \subset U_{{\cal F}_{{\cal E}_3} {\cal F}_{%
{\cal E}_2}}(sp(2)) \subset U_{{\cal F}_{{\cal E}_3}{\cal F}_{{\cal E}_2}%
{\cal F}_{{\cal E}_1}}(sp(3)).
$$
Still it is based on the set of symplectic regular subalgebras. Thus the
chain (\ref{f-chain}) is {\sl proper} for the algebra $sp(3)$ (in contrast
with the improper chains described in \cite{KLO}). The carrier subalgebra
for (\ref{f-chain}) is equivalent to $B^+(sp(3))$. Thus this chain of
extended twists is {\sl full}.

\section{Conclusions}

The importance of explicitly defined twisting elements is that the
corresponding quantum $R$-matrices can be calculated in any representation.
The general formula for the triangular universal $R $-matrix in the case
of $sp(3)$ algebra is as
follows:
\begin{equation}
\label{rmat}
{\cal R}=({\cal F}_{{\cal G}_{0 \prec 2}})_{21} ({\cal F}_{{\cal G}_{0 \prec
2}})^{-1}.
\end{equation}

The classical limit can be obtained after the rescaling (and turning
$\xi$ to zero)
$$
E_{i+k}\longrightarrow \xi\eta_{k} E_{i+k};\quad k\ge i.
$$

And the classical $r$-matrix has the form
$$
r=\sum_{k=1}^{3}\eta_{k}(H_{kk}\land E_{k+k}+ \kappa_{k}\sum_{i=k+1}^{3}
E_{k-i}\land E_{k+i}),
$$
where $\kappa_{k}$ (the discrete parameter equal to zero or one) indicates
whether we perform in the subalgebra the extended jordanian twist or the
pure jordanian twist. Notice, that only two $\kappa _{k}$'s can appear in
this formula ($\kappa _{1}$ and $\kappa _{2}$). And in the case of a full
chain we must keep $\kappa _{k}=1 $, for any $k $. Applying the operator
$ \exp (\alpha {\rm ad}(H_{ss}\otimes1 + 1\otimes H_{ss}))$ to quantum
$R$-matrix (\ref{rmat}) we can change the values of
$\eta_{s}$ $(\eta_s \rightarrow \eta_{s}(1+\alpha))$.

These results can be generalized for the case of an arbitrary symplectic
simple algebra $sp(N)$. On any step of a chain the transformations of
basis
leading to the deformed primitive carrier subspaces can be written
explicitly. The corresponding twisting element has the following form:
\begin{equation}
\label{f-genchain}
\begin{array}{rcl}
{\cal F}_{{\cal G}_{0 \prec (N-1)}} & = & e^{H_{NN}\otimes
\sigma^{(N-1)}_{N+N}}\cdot \\
&  & \cdot e^{E^{(N-2)}_{(N-1)-N}\otimes E^{(N-2)}_{(N-1)+N} e^{-
\frac{1}{2}\sigma^{(N-2)}_{(N-1)+(N-1)}}} e^{H_{(N-1)(N-1)}\otimes
\sigma^{(N-2)}_{(N-1)+(N-1)}}\cdot \\  &  & \vdots \\
&  & \cdot e^{E_{2-N}\otimes E_{2+N}e^{-\frac 12\sigma^{\prime}_{2+2}}}\dots
e^{E_{2-3}\otimes E_{2+3}e^{-\frac 12\sigma^{\prime}_{2+2}}}
e^{H_{22}\otimes \sigma^{\prime}_{2+2}}\cdot \\
&  & \cdot e^{E_{1-N}\otimes E_{1+N}e^{-\frac 12\sigma _{1+1}}}\dots
e^{E_{1-2}\otimes E_{1+2}e^{-\frac 12\sigma_{1+1}}} e^{H_{11}\otimes \sigma
_{1+1}} \\
&  &
\end{array}
\end{equation}
Here the elements $\sigma _{K+K}^{(L)}$ are the analogs of $\sigma
_{2+2}^{\prime }$ and $\sigma _{3+3}^{\prime \prime }$ introduced in 
Section 2.

Thus we can operate with the explicit expression for the corresponding
quantum $R$-matrix
\begin{equation}
\label{qr-gen}{\cal R}=({\cal F}_{{\cal G}_{0\prec (N-1)}})_{21}({\cal F}_{%
{\cal G}_{0\prec (N-1)}})^{-1}.
\end{equation}
According to the general properties of chains \cite{KLO} every step in (\ref
{f-genchain}) can bare its own continuous deformation parameter. Thus the $R$%
-matrix has a natural multiparametric form. To get the classical $r$-matrix
one must fix the quantization curve that leads to the classical universal
enveloping algebra $U(sp(N))$. The general position of the classical limit
curve is obtained when all the parameters are proportional to, say, the
length of the curve. In this case the most general classical $r$-matrix
(whose quantum analog is (\ref{qr-gen})) can be written in the following
form:
$$
r=\sum_{k=1}^N\eta _k(H_{kk}\land E_{k+k}+\kappa
_k\sum_{i=k+1}^NE_{k-i}\land E_{k+i}).
$$
The continuous parameters $\eta _k$ describe the ratios of the
canonical parameters while the discrete parameters $\kappa _k$ indicate the
possibility to switch off the extension factors in the chain.

\end{document}